\documentclass{amsart}

\usepackage{amsmath}
\usepackage{amssymb}

\def\pa{\partial}
\def\var{\varepsilon}
\def\R{{\Bbb R}}

\newtheorem{Thm}{Theorem}
\newtheorem{Prop}[Thm]{Proposition}

\newtheorem{Lem}{Lemma}
\newtheorem{Rmk}{Remark}

\def\dd{\varepsilon}
\def\g{\gamma}

\begin{document}

\title[Entropy dissipation estimates for the Landau equation]{Some remarks about the Link between the Fisher information and Landau or Landau-Fermi-Dirac entropy dissipation }

\author{L. Desvillettes}

\address{Universit\'e Paris Cit\'e and Sorbonne Universit\'e, CNRS, IUF, Institut de Math\'ematiques de Jussieu-Paris Rive Gauche (desvillettes@math.univ-paris-diderot.fr)}

\begin{abstract}
We present in this work variants of existing estimates for the Landau or Landau-Fermi-Dirac entropy dissipation, in terms of Fisher information, in the hard potential case. The specificity of those variants is that the entropy is never used in the estimates (in order to control possible concentrations on a zero measure set). The proofs are significantly simplified with respect to previous papers on the subject
\end{abstract}


\maketitle

\section{Introduction
} \label{intro}

\subsection{Description of the Landau operator and its entropy dissipation}


We write down the Landau operator appearing in kinetic theory (cf. \cite{chapman,lifschitz}),
defined by
\begin{equation} \label{operator}
Q_{\gamma}(f,f)(v) = \nabla_{v} \cdot   \left \{ \int_{{\Bbb R}^3} \,
 \, |v-w|^{2+\gamma}\,  \Pi (v-w) \,
\bigg( f(w) \,\nabla {f}(v) - f(v) \,
\nabla {f}(w) \bigg)\,  dw\right \} ,
\end{equation}
where 
\begin{equation} \label{pii}
\Pi_{ij}(z) = \delta_{ij} - \frac{z_i z_j}{|z|^2} 
\end{equation}
is the $i,j$-component of the orthogonal projection $\Pi$
 onto $z^\bot := \{ y\, / \; y\cdot z =0 \}.$
\medskip

In this paper, we focus on the case of so-called 
 Maxwell molecules, hard spheres or hard potentials, that is
\begin{equation}\label{hshp}
\gamma \in [0,1] . 
\end{equation}
\medskip


 
The Landau operator satisfies (at the formal level) the conservation of 
mass, momentum and kinetic energy, that is, for any $f := f(v)\ge 0$:
\begin{equation} \label{conserv}
 \int_{\R^3} Q_{\gamma}(f,f)(v) \, \left(\begin{array}{c} 1\\ v_i\\ |v|^2/2 \end{array} \right)\, dv =  \left(\begin{array}{c} 0\\ 0 \\ 0 \end{array} \right). 
\end{equation}
We also have the following property of nonnegativity of the so-called entropy dissipation associated to the Landau operator, 
namely (at the formal level), for any $f:=f(v) >0$,
\begin{equation}\label{mded} 
D_{\gamma}(f)  := -\int_{{\Bbb R}^3}  Q_{\gamma}(f,f)(v) \ln f(v) \, dv 
\end{equation}
$$=\frac12\, \iint_{{\Bbb R}^3\times {\Bbb R}^3} 
f(v)\,f(w)\,
|v-w|^{2 + \gamma} \, \Pi(v-w)\, \left ( \frac{\nabla f}{f}(v) -
\frac{\nabla f}{f}(w) \right ) $$
 $$\cdot \,   \left( \frac{\nabla f}{f}(v) -
\frac{\nabla f}{f} (w) \right )\, dv dw \ge 0. $$
\medskip

In the last years, several works (cf. \cite{DV2}, \cite{FPL}, \cite{braga1}, \cite{X1}, \cite{X2}) have been devoted to the proof of estimates from below of the quantity $D_{\gamma}(f)$ in terms of 
quantities like the Fisher information $F(f) := \int_{\R^3} \frac{|\nabla f(v)|^2}{f(v)}\, dv$, or the relative (with respect to a Maxwellian function $M:=M(v)$) Fisher information $F_{rel}(f) := \int_{\R^3} \bigg| \frac{\nabla f}{f}(v) - \frac{\nabla M}{M}(v) \bigg|^2 \, f(v) \,dv$.
\medskip

Those estimates (typically the ones related to $F(f)$) yield 
 results of regularity for the solutions to the Landau equation (cf. \cite{FPL}, \cite{GGIV}, \cite{DHJ}). When they are complemented with the logarithmic Sobolev inequality of Gross (cf. \cite{gross}), they also 
 (typically the ones related to $F_{rel}(f)$) yield results for the large time behavior of the Landau equation (cf. \cite{DV2}, \cite{CDH}), following 
the lines of the entropy-entropy dissipation method introduced in kinetic theory in the 90s by \cite{carcar1}, \cite{carcar2}, \cite{toscani} (cf. \cite{DMV} for more historical details and an explanation of the general context). 
\medskip

We write down a typical estimate, taken from Corollary 2.5 of \cite{braga1},  which holds when $\g$ satisfies (\ref{hshp}), and $f > 0$, 
$$  \int_{\R^3} f(v) \, \left( \begin{array}{c} 1 \\ v \\ |v|^2 \end{array}  \right) \, dv = \left( \begin{array}{c} 1 \\ 0 \\ 3 \end{array}  \right)$$
 (and when the terms appearing in the inequality are defined):
\begin{equation} \label{abstract6}
\int_{\R^3} \frac{|\nabla f (v)|^2}{f(v)} \, dv \le 3072\,\Delta (f)^{-2}\, 
\bigg\{8448 + 48 \, \sqrt{1+\pi} \, D_{\gamma}(f) \, ||f||_{L^2(\R^3)} \bigg\}, 
\end{equation}
where
\begin{equation} \label{dea}
\Delta (f) := {\hbox{ Det}} \, \bigg(  \int_{\R^3}  f(w)\, (1+ |w|^2)^{-1/2}\, 
 \left[ \begin{array}{ccc}
 1 & w_i & w_j \\
 w_i & w_i^2 &  w_i\,w_j \\
 w_j & w_i\,w_j &  w_j^2 \end{array}\, \right] \,dw\, \bigg) . 
\end{equation}

 The purpose of this paper is to show that it is possible to obtain estimates like (\ref{abstract6}) with the following improvements:
\begin{itemize}
\item
No terms involving a determinant like (\ref{dea}) appears;
\item
The numerical constants are much smaller than the ones appearing in (\ref{abstract6});
\item
The proof is short;
\item
These new estimates extend somewhat  to the case of the Landau-Fermi-Dirac entropy dissipation.
\end{itemize}

 The price to pay in order to get those improvements is the following restriction: the estimates that we will write down hold only when $D_{\gamma}(f)$ is not too large. We will see however that this is not a problem when one wishes to prove the exponential decay towards equilibrium with explicit rate
for the solutions to the Landau equation with hard potentials or hard spheres.  
\medskip

 We notice that the statement and computations presented below have the same flavour as the material presented in chapter 5 of \cite{DV2} (more generally, the methodology for getting estimates relating $D_{\gamma}$ and $F$ or $F_{rel}$ was introduced in \cite{DV2}, recalling some of the computations of \cite{desv89}). The conditions on the cross section are however more general in the present paper, and the proofs are simplified a lot with respect to what is presented in \cite{DV2}.  
\medskip

Our main result writes:

\begin{Thm} \label{main}
We consider $\gamma$ satisfying (\ref{hshp}). Then 
for all $f:= f(v) > 0$ (sufficiently smooth for the terms in the estimate to make sense) satisfying 
\begin{equation}\label{norm}
 \int_{\R^3} f(v) \, \left( \begin{array}{c} 1 \\ v \\ |v|^2 \end{array}  \right) \, dv = \left( \begin{array}{c} 1 \\ 0 \\ 3 \end{array}  \right) ,
\end{equation}
and
such that
$$ \bigg( \int_{\R^3} |v|^4 \, f(v) \, dv \bigg)  \,\left( 1 + 2\sqrt{\pi}\, ||f||_{L^2(\R^3)} \right) \,  D_{\gamma}(f) \le \frac{27}{32} , $$
 the following estimate holds: 
$$ \int_{\R^3} \bigg| \frac{\nabla f (v)}{f(v)} + v\, \bigg|^2 \, f(v)\, dv \le\bigg[  \frac{64}{9} \, \bigg(  \int_{\R^3} |v|^4 \, f(v) \, dv \bigg)  \,\left( 1 + 2\sqrt{\pi}\, ||f||_{L^2(\R^3)} \right) $$
\begin{equation}\label{resth}
 +\, 48 + 32 \sqrt{\pi} \, ||\, |\cdot|^2\,f||_{L^2(\R^3)}  \bigg] \, D_{\gamma}(f) .
\end{equation}

\end{Thm}

\begin{Rmk}
Inequality (\ref{resth}) can be further simplified, and be transformed for example in
$$  \int_{\R^3} \bigg| \frac{\nabla f (v)}{f(v)} + v \bigg|^2 \, f(v)\,dv  \le 200\,  ||f||_{L^2_6(\R^3)}^2 \, D_{\gamma}(f) ,$$ 
as soon as 
$$  ||f||_{L^2_6(\R^3)}^2 \, D_{\gamma}(f) \le 0.062, $$
where we denote from now on (for any $p \in [1, \infty[$ and $q\in\R$)
$$ || f ||_{L^p_q(\R^3)}^p := \int_{\R^3} (1 + |v|^2)^{pq/2} \, |f(v)|^p \, dv . $$ 
Indeed, one can check thanks to H\"older's inequality that (when $f$ satisfies (\ref{norm}))
$$ \bigg( \int_{\R^3} |v|^4 \, f(v) \, dv \bigg)  \,\left( 1 + 2\sqrt{\pi}\, ||f||_{L^2(\R^3)} \right) \le \left(2\pi\,\sqrt{\pi} + \frac{\pi^2}4 \right)\, 
 ||f||_{L^2_6(\R^3)}^2. $$
We see at this point that the numerical constant (that is, 200) is indeed much smaller than the ones appearing in inequality (\ref{abstract6}).
\end{Rmk}
\medskip

 Section \ref{sec2} is devoted to the proof of Theorem \ref{main}. Then, we show in section~\ref{sec3} how this result can be used to obtain in a quick way the  explicit exponential decay towards equilibrium for the solutions to the Landau equation with hard potentials or hard spheres. Finally, we show in section \ref{sec4} that Theorem \ref{main} can be extended somewhat in the case of the Landau-Fermi-Dirac equation. 

\section{Proof of Theorem \ref{main}}\label{sec2}

We start here the
\medskip

\noindent
{\bf{Proof of Theorem \ref{main}}} :
We denote in this proof,  for $i =1,2,3$,
\begin{equation}\label{ii} 
 I_i(f) := \int_{\R^3} f(v)\, v_i^2\, dv, 
\end{equation}
so that $I_1(f) + I_2(f) + I_3(f)=3$. 
\medskip

Up to the use of a linear isometry of $\R^3$, we assume like in \cite{DV2} without loss of generality that, for $i,j =1,2,3$ such that $i \neq j$,
\begin{equation}\label{rot} 
\int_{\R^3} f(v)\, v_i\,v_j\,  dv = 0 . 
\end{equation}
Note that the Landau operator, like the Boltzmann operator (from which it is a singular limit: the so called grazing collisions limit) is indeed built to be Galilean-invariant. 

 We first observe that 
$$ D_{\gamma}(f) = \frac12 \sum_{i=1}^3  \int\int_{\R^3 \times \R^3}  f(v)\,f(w)\, |v-w|^{\gamma}\, [q_{i, i+1}^f(v,w)]^2 \, dv dw , $$
where for $i,j =1,2,3$,
 \begin{equation}\label{qij}  
q_{i,j}^f(v,w) :=  (v_i - w_i) \, \left( \frac{\pa_j f}f (v) - \frac{\pa_j f}f (w) \right) -
    (v_j - w_j) \, \left( \frac{\pa_i f}f (v) - \frac{\pa_i f}f (w) \right), 
\end{equation}
and  by definition $q_{i,4}^f := q_{i,1}^f$.
\medskip

Using the normalizations (\ref{norm}) and (\ref{rot}), we multiply (\ref{qij}) by $f(w)$ and $w_i\, f(w)$ and integrate on $\R^3$ in the $w$-variable.
 We get for $i \neq j$ the identities:
\begin{equation}\label{prin1}  
 v_i\, \frac{\pa_j f}f (v) - v_j \,  \frac{\pa_i f}f (v) = \int_{\R^3}  f(w)\, q_{i,j}^f(v,w) \, dw, 
\end{equation}
\begin{equation}\label{prin2} 
  I_i(f)\,  \frac{\pa_j f}f (v) + v_j  =  -  \int_{\R^3} w_i\, f(w)\, q_{i,j}^f(v,w) \, dw . 
\end{equation}
It is convenient to write (\ref{prin2}) under the form
\begin{equation}\label{prin2bis} 
  I_i(f)\,  \left( \frac{\pa_j f}f (v) + v_j  \right) = (I_i(f) - 1)\, v_j -  \int_{\R^3} w_i\, f(w)\, q_{i,j}^f(v,w) \, dw . 
\end{equation}
  As a consequence, using the Cauchy-Schwarz inequality,
$$ I_i(f)^2 \int_{\R^3} \bigg|  \frac{\pa_j f}f (v) + v_j  \bigg|^2\, f(v) \, dv  \le 2 \, (I_i(f) - 1)^2\,I_j(f) $$
$$+\, 2  \int_{\R^3}  f(v) \, \bigg|  \int_{\R^3} w_i\, f(w)\, q_{i,j}^f(v,w) \, dw  \bigg|^2 \, dv $$
 $$ \le  2 \, (I_i(f) - 1)^2\,I_j(f) + \, 2   \int_{\R^3}  f(v)  \,\bigg[ \int_{\R^3}  f(w)\, |v-w|^{\gamma}\, |q_{i,j}^f(v,w)|^2 \, dw \bigg] $$
$$ \times\,  \,\bigg[ \int_{\R^3}  f(w)\, |v-w|^{- \gamma} \, w_i^2\, dw \bigg] \, dv $$
\begin{equation}\label{uti} 
 \le  2 \, (I_i(f) - 1)^2\,I_j(f) +  \, 2   \int_{\R^3}  f(v)  \,\bigg[ \int_{\R^3}  f(w)\, |v-w|^{\gamma}\, |q_{i,j}^f(v,w)|^2 \, dw \bigg] \, dv\, \sigma_{\gamma}(f) , 
\end{equation}
where 
\begin{equation}\label{sf}
 \sigma_{\gamma}(f) := \sup_{v \in \R^3} \int_{\R^3}  f(w)\, |v-w|^{-\gamma} \, |w|^2\, dw . 
\end{equation}
Multiplying identity (\ref{prin1}) by $v_i\,v_j\,f(v)$  and integrating on $\R^3$ in the $v$-variable,
we see that
$$ I_i(f) - I_j(f) =  \int\int_{\R^3 \times \R^3}  f(v)\,f(w)\, v_i\, v_j \, q_{i,j}^f(v,w) \, dv dw $$
so that, using again the Cauchy-Schwarz inequality, 
$$ | I_j(f) - I_i(f)| \le  \bigg[\int \int_{\R^3 \times \R^3} f(v)\,  f(w)\, |v-w|^{\gamma}\, |q_{i,j}^f(v,w)|^2 \, dv dw \bigg]^{1/2}$$
$$ 
\times  \,\bigg[ \int\int_{\R^3 \times \R^3}  f(v)\, f(w)\, |v-w|^{- \gamma} \, v_i^2\, v_j^2 \, dv dw \bigg]^{1/2} $$
$$ \le \sqrt{\frac23} \, D_{\gamma}(f)^{1/2} \, s_{\gamma}(f)^{1/2}, $$
 where 
\begin{equation}\label{vsf}
  s_{\gamma}(f) := \int\int_{\R^3 \times \R^3}  f(v)\, f(w)\, |v-w|^{-\gamma} \, |v|^4\, dv dw . 
\end{equation}
As a consequence,
$$ |I_i(f) - 1|^2 = \left(\frac23 I_i(f) - \frac13 I_{i+1}(f) - \frac13 I_{i+2}(f) \right)^2 $$
$$ \le \frac29 (I_i(f) -  I_{i+1}(f))^2 +  \frac29 (I_i(f) -  I_{i+2}(f))^2 \le \frac8{27}\,  D_{\gamma}(f) \,  s_{\gamma}(f). $$
We deduce from this estimate that
$ I_i(f) \ge \frac12$ as soon as 
\begin{equation}\label{conddp}
  D_{\gamma}(f) \,  s_{\gamma}(f) \le \frac{27}{32} ,
\end{equation}
since (\ref{conddp}) implies that $ |I_i(f) - 1| \le \frac12$. 
\par
Under condition (\ref{conddp}), we see therefore that (\ref{uti}) leads to
$$  \int_{\R^3} \bigg|  \frac{\pa_j f}f (v) + v_j  \bigg|^2\, f(v) \, dv
  \le 8\, (I_i(f) - 1)^2\,I_j(f) $$ 
$$ +\,  8   \int_{\R^3}  f(v)  \,\bigg[ \int_{\R^3}  f(w)\, |v-w|^{\g}\, |q_{i,j}^f(v,w)|^2 \, dw \bigg] \, \sigma_{\g}(f) ,$$
so that finally, summing from $j= 1$ to $3$ (and choosing $i=j+1$)
 \begin{equation}\label{aa1}
 \int_{\R^3} \bigg|  \frac{\nabla f}f (v) + v  \bigg|^2\, f(v) \, dv
   \le  \left(\frac{64}{9} \, s_{\g}(f) + 16 \, \sigma_{\g}(f) \right)\,  D_{\g}(f) . 
\end{equation}
Then we observe that 
$$ s_{\gamma}(f) \le \int\int_{|v-w| \ge 1} |v|^4 \, f(v)\, f(w) \, dv dw + \int\int_{|v-w| \le 1} |v|^4 \, f(v)\, f(w) \, |v-w|^{- \gamma} \, dv dw $$
$$ \le \bigg( \int_{\R^3} |v|^4 \, f(v) \, dv \bigg)  \, \left(1 + \sup_{v \in \R^3} \int_{|v-w| \le 1}  f(w)\,  |v-w|^{- \gamma} \,  dw \right) $$
 \begin{equation}\label{aa2}
 \le\bigg(  \int_{\R^3} |v|^4 \, f(v) \, dv \bigg)   \,\left( 1 + 2\sqrt{\pi}\, ||f||_{L^2(\R^3)} \right) , 
\end{equation}
and that
$$ \sigma_{\g}(f)  \le \sup_{v\in\R^3}   \int_{|v-w| \ge 1} |w|^2 \, f(w) \, dw
 +  \sup_{v\in\R^3}   \int_{|v-w| \le 1} |w|^2 \,  |v-w|^{- \gamma} \, f(w) \, dw $$
\begin{equation}\label{aa3}
 \le 3 + 2\sqrt{\pi} \, || \,  |\cdot|^2\, f||_{L^2(\R^3)} . 
\end{equation}
We conclude the proof of Theorem \ref{main} by using estimates (\ref{aa2}) and (\ref{aa3}) in estimate (\ref{aa1}). 
$\square$

\section{Application to the large time behavior of the Landau equation} \label{sec3}

 We show in this section that Theorem \ref{main} quickly leads
 to a result of exponential convergence towards equilibrium with explicit rate, for
the Landau equation in the case of hard potentials or hard spheres. 
\medskip

Note that in \cite{DV2}, such a result is proven in the so-called overMaxwelian case, which includes the case of
 Maxwell molecules, but not the case of hard potentials or hard spheres, because in those cases, the cross section cancels at some point. Those cases are treated as a special case (when the quantum parameter $\var$ tends to $0$) of the results obtained in \cite{ABDL} for the Landau-Fermi-Dirac equation. What we now present is a new (somewhat simpler) proof.
\medskip

 We start with a simple Lemma which replaces here Proposition 6 and its consequences in \cite{DV2}:

\begin{Lem} \label{l1}
Suppose that $H \in C^1(\R_+; \R_+)$ and $D \in C(\R_+; \R_+)$ are such that 
$$ - H' = D . $$
 We also assume that for some real numbers $q,c_0 >0$ 
 \begin{equation}\label{qc}
 D \le q \qquad \Rightarrow \qquad D \ge c_0\, H.
\end{equation}
Then,  for all $t \ge H(0)/q$, the following estimate holds:
 \begin{equation}\label{ll}
  H(t) \le H(0) \, \exp( c_0 \, H(0)/q) \, e^{- c_0\,t} . 
\end{equation}
\end{Lem}

In the sequel, $H$ and $D$ will be the respective relative entropy and entropy dissipation of the Landau
and  the Landau-Fermi-Dirac equations.
\medskip

\noindent
{\bf{Proof of Lemma \ref{l1}}} : 
We first observe that if for some ${\bar{t}} \ge 0$, one has $H(\bar{t}) \le q/c_0$, then for all $t \ge {\bar{t}}$,
 \begin{equation}\label{ll2}
  H(t) \le q/c_0
\end{equation}
since $H$ is nonincreasing. 
\par
Then either $D(t) \le q$ and thanks to (\ref{qc}), $D(t) \ge c_0\, H(t)$; or $D(t) \ge q \ge c_0\, H(t)$ thanks to (\ref{ll2}).
\medskip

In both cases, $ - H'(t) \ge c_0\,H(t)$ for $t \ge \bar{t}$.
Consequently, 
\begin{equation}\label{ll3}
 \forall t \ge \bar{t}, \qquad  H(t) \le  H(\bar{t}) \, e^{- c_0\, (t - \bar{t})} .
\end{equation}
We apply estimate (\ref{ll3}) for 
$$ \bar{t} :=  \inf\{ \tau \ge 0, H(\tau) \le q/c_0 \}. $$
We see that $H(s) \ge q/c_0$ for $s \in ]0, \bar{t}]$ (this interval is empty if $\bar{t}=0$, and might be of infinite length). 
Then, for   $s \in ]0, \bar{t}]$,  if $D(s) \le q$ we see that $D(s) \ge c_0\, H(s) \ge q$ (because of (\ref{qc})). Therefore $D(s) \ge q$ on $[0, \bar{t}]$, so that  $H'(s) \le - q$ on the same interval, and finally 
$$ 0 \le H(\bar{t}) \le H(0) - q\, \bar{t} , $$
so that $\bar{t} \le H(0)/q$.
\medskip

Coming back to (\ref{ll3}), we see that for $t \ge H(0)/q \, \, (\ge \bar{t})$,
$$ H(t) \le H(\bar{t}) \, \exp(c_0\, H(0)/q)\, e^{- c_0\,t}. $$
We conclude that estimate (\ref{ll}) holds (using the fact that $H$ is a nonincreasing function).
$\square$
\bigskip

We now explain how this Lemma (together with Theorem \ref{main}) yields a result of exponential convergence towards equilibrium with explicit rate for
the Landau equation.
\medskip
 
Using Theorem 5 of \cite{DV1}, we see that when $\gamma \in ]0,1]$ (that is, in the case of hard potentials or hard spheres), and when the initial datum $f_{in} \ge 0$ lies in
$L^2_3(\R^3) \cap L^1_3(\R^3)$ (beware that the definition of $L^p_q$ is different here and in \cite{DV1}), then there exists a solution $f:=f(t,v)\ge 0$ to the (spatially homogeneous) Landau equation 
$\pa_t f = Q_{\g}(f, f)$ (with initial datum $f_{in}$) lying in  $L^{\infty}([1, \infty[ ; L^2_6(\R^3))$, and
satisfying the entropy identity in the strong sense (this is actually true if the initial datum is bounded below by a Maxwellian and has a sufficient number of moments; if the initial datum does not fulfill those requirements, one has to approximate it by initial data fulfilling it): 
$$  - \frac{d}{dt} H(f(t, \cdot)| M) =  D_{\gamma}(f(t, \cdot)) ,$$
where 
$$ H(f|M) := \int_{\R^3} f(v)\, \log f(v) dv - \int_{\R^3} M(v)\, \log M(v) dv, $$
and $M (v) = (2\pi)^{-3/2}\, e^{- \frac{|v|^2}2}$.
 \medskip

At this point we assume the normalization $\int_{\R^3} f_{in}(v) \, \left( \begin{array}{c} 1 \\ v \\ |v|^2 \end{array}  \right) \, dv = \left( \begin{array}{c} 1 \\ 0 \\ 3 \end{array}  \right)$,  and apply Theorem \ref{main} (more precisely Remark 1) with  $f: = f(t, \cdot)$ (the considered solution to the Landau equation).
\medskip

 It is then possible to use Lemma \ref{l1} (with a time translation), using $$H: = H(f(t, \cdot)| M), \qquad D:= D_{\gamma}(f(t, \cdot)), $$
  $$ q:= 0.22\, ||f||_{ L^{\infty}([1, \infty[ ; L^2_6(\R^3))}^{-2}, \qquad c_0 : = c_1\, 200^{-1}\,  ||f||_{ L^{\infty}([1, \infty[ ; L^2_6(\R^3))}^{-2}, $$
where $c_1$ is the constant in the logarithmic Sobolev inequality of Gross \cite{gross}. 
\medskip

We obtain thanks to Lemma \ref{l1}
 that this solution decays exponentially fast, with an explicit rate, towards the equilibrium $M$.
\medskip

It is possible to obtain the same result without assuming the normalization (but with a different Maxwellian equilibrium $M$) thanks to the 
change of unknown $f \to a\,f( b\, (\cdot - u))$, where $a,b>0$, and $u \in \R^3$.
\medskip

\begin{Rmk}
A similar result can be obtained for the case of Maxwell molecules ($\g =0$), but in that case, direct explicit computations can be exploited in order
to get the exact rate of exponential convergence. 
\par
The case of soft potentials $\g \in [-3,0[$ is somewhat different. Moments estimates still exist (cf. \cite{FPL} and \cite{CDH}), as well as smoothness estimates when $\g \in [-2, 0[$ (case of moderately soft potentials, cf. for example \cite{kcw} and \cite{FPL}), but they are not uniform in time. It is nevertheless possible (at least in the Coulomb case $\g = -3$) to obtain a stretched exponential estimate of convergence towards equilibrium (cf. \cite{CDH}), using methods related to those of \cite{toscani_villani}. However, it is not obvious to know if such a result can be recovered with a simplified proof using
arguments similar to those used in this paper. 
\end{Rmk}

\section{Extension to the Landau-Fermi-Dirac equation} \label{sec4}

The Landau Fermi Dirac operator (cf. \cite{bag}, \cite{ABL}, \cite{ABDL}, \cite{ABDL2}) is 
defined by

$$Q_{\g, LFD}(f,f)(v) = \nabla_{v} \cdot   \left \{ \int_{{\Bbb R}^3} \, |v-w|^{2+\g} \, \Pi(v-w)\,
\bigg( f(w) \, (1 - \var\, f(w))\,\nabla {f}(v) \right. $$
\begin{equation} \label{operatorLFD}
\left.  - f(v) \, (1 - \var\, f(v))
\nabla {f}(w) \bigg)\,  dw\right \} ,
\end{equation}
where $\Pi$ is defined by (\ref{pii}), and $\var>0$ is a parameter which mesures 
the amount of quantum effects (due to the Pauli exclusion principle). 
\medskip

For this operator, we still have (at the formal level) the conservation of 
mass, momentum and kinetic energy, that is, for any $f := f(v)\ge 0$:
\begin{equation} \label{conserv_lfd}
 \int_{\R^3} Q_{\g, LFD} (f,f)(v) \, \left(\begin{array}{c} 1\\ v_i\\ |v|^2/2 \end{array} \right)\, dv = \left(\begin{array}{c} 0 \\ 0 \\ 0 \end{array} \right).
\end{equation}
The  entropy dissipation associated to the Landau-Fermi-Dirac operator also remains a nonnegative quantity, but it is not identical to the entropy dissipation associated to the Landau operator. Note that the entropy itself is different in the 
classical and quantum case (cf. (\ref{aaa})). More precisely (at the formal level), for any $f:=f(v) >0$, it writes
\begin{equation}\label{mded} 
D_{\g, LFD}(f)  := -\int_{{\Bbb R}^3}  Q_{\g, LFD}(f,f)(v)\, \left[ \ln f(v) - \ln (1 - \var\,f(v)) \right] \, dv 
\end{equation}
$$=\frac12\, \iint_{{\Bbb R}^3\times {\Bbb R}^3} 
f(v)\, \, (1 - \var\, f(v)) \,f(w)\, \, (1 - \var\, f(w)) |v-w|^{2+\g} \, \Pi(v-w) $$ 
$$ \left ( \frac{\nabla f}{f\,(1 - \var\,f)}(v) -
\frac{\nabla f}{f\,(1 - \var \,f)}(w) \right ) \, \cdot\, \left( \frac{\nabla f}{f\, (1 - \var\,f)}(v) - \frac{\nabla f}{f\,(1 - \var f)} (w) \right )\, dv dw \ge 0. $$
\medskip

 In \cite{ABDL}, it is shown that this quantity can be related to a quantity close to the Fisher information (more precisely, a quantum variant of $F$ or $F_{rel}$). More precisely, the following estimate is presented (Proposition 2.12 of \cite{ABDL}): 
for $\g$ satisfying assumption (\ref{hshp}) and for $f : = f(v) >0$ such that (\ref{norm}) holds and 
\begin{equation}\label{ka}
 \forall v \in \R^3   \qquad   1 - \var\, f(v)\ge \kappa_0, 
\end{equation}
one has (when the quantities below are well defined)
$$ \int_{\R^{3}}\left|\frac{\nabla f(v)}{f(v)(1-\dd f(v))}-Kv\right|^{2} \, f(v)\, d v \leq 510  \,\,(\min_{i}I_{i}(f))^{-3}\, {\kappa_{0}^{-2}}\,\max(1, B_{\g}(f))$$ 
 \begin{equation}\label{eq:entrop} 
 \times \, \max\left(1,\int f(v)\, (1+ |v|^2)^{1+\frac{\g}{2} }\,dv \right) \, {J}_{\gamma}(f)\, {D}_{\g, LFD}(f),
\end{equation}
where  $K := \dfrac{1}{\dd}\displaystyle\int_{\R^{3}}\ln(1-\dd f(v) ) \, dv$, the directional temperature  $I_i(f)$ is defined in (\ref{ii}), 
$${B}_{\g}(f) := \bigg(\min_{i\neq j}\inf_{\sigma \in S^{1}}\int_{\R^{3}}\left|\sigma_{1}\frac{v_{i}}{\langle v\rangle}-\sigma_{2}\frac{v_{j}}{\langle v\rangle}\right|^{2} f(v) \, dv \bigg)^{-1}, $$
and 
$${J}_{\gamma}(f) :=\sup_{v \in \R^{3}}\langle v\rangle^{\gamma}\int_{\R^{3}}f(w)\, |w-v|^{-\gamma}\, (1+ |w|^2) \, d w .$$
\medskip

In the same spirit as in Section \ref{sec2}, we propose an estimate which is close to (\ref{mded}), but somewhat simpler, and much easier to prove. It writes:

\begin{Prop} \label{main2}
We consider   $\gamma$ satisfying (\ref{hshp}). Then
for all $f:= f(v) > 0$ (sufficiently smooth for the terms in the estimate to make sense) satisfying the normalization (\ref{norm}), such that 
\begin{equation}\label{kap}
 \forall v \in \R^3   \qquad   1 - \var\, f(v)\ge \kappa_0, 
\end{equation}
and such that
$$ \bigg( \int_{\R^3} |v|^4 \, f(v) \, dv \bigg)  \,\left( 1 + 2\sqrt{\pi}\, ||f||_{L^2(\R^3)} \right) \, \kappa_0^{-1}\,  D_{\g, LFD}(f) \le \frac{27}{32} , $$
 the following estimate holds: 
$$ \int_{\R^3} \bigg| \frac{\nabla f (v)}{f(v) \,(1-\dd f(v))} - K\, v \bigg|^2 \, f(v)\, dv \le 
 \bigg[  \frac{32}3 \, \kappa_0^{-1}\, \bigg( \int_{\R^3} |v|^4 \, f(v) \, dv \bigg)  \,\left( 1 + 2\sqrt{\pi}\, ||f||_{L^2(\R^3)} \right) $$
\begin{equation}\label{rlfd_LFD}
 +\, 216\, \left( 1 + 2\sqrt{\pi}\, ||f||_{L^2(\R^3)} \right)  + 24\, (3+ 2\sqrt{\pi}\, ||\, |\cdot|^2\, f||_{L^2(\R^3)} )\, \bigg] \, \kappa_0^{-2}\, D_{\g, LFD}(f) ,
\end{equation}
where we recall that  $K := \dfrac{1}{\dd}\displaystyle\int_{\R^{3}}\ln(1-\dd f(v) ) \, dv$.
\end{Prop}

\begin{Rmk}
An important feature of the above estimate is that it is uniform with respect to $\var$ (for $\var >0$ such that (\ref{kap}) holds, like all the estimates provided in \cite{ABDL}). One  can check that when $\var \to 0$, Proposition \ref{main2} becomes Theorem \ref{main}, up to the numerical constants. Remark 1 can also  be adapted to the case of the Landau-Fermi-Dirac operator. 
\end{Rmk}
\medskip

\noindent
{\bf{Proof of Proposition \ref{main2}}} :
 We first observe that 
$$ D_{\g, LFD}(f) = \frac12 \sum_{i=1}^3  \int\int_{\R^3 \times \R^3}  f(v)\,(1-\dd f(v)) \,f(w)\, (1-\dd f(w))$$
$$ \times\,  |v-w|^{\gamma} \, [r_{i, i+1}^f(v,w)]^2 \, dv dw , $$
where for $i,j =1,2,3$,
$$ r_{i,j}^f(v,w) :=  (v_i - w_i) \, \left( \frac{\pa_j f}{f \,(1-\dd f)}(v) - \frac{\pa_j f}{f\,(1-\dd f)}(w) \right)$$
\begin{equation}\label{qijlfd}   
 -    (v_j - w_j) \, \left( \frac{\pa_i f}{f\,(1-\dd f)}(v) - \frac{\pa_i f}{f\,(1-\dd f)}(w) \right), 
\end{equation}
and,  by definition, $r_{i,4}^f := r_{i,1}^f$.
\medskip

Using the normalization (\ref{norm}) and assuming without loss of generality that (\ref{rot}) holds,
 we multiply (\ref{qijlfd}) by $f(w)$ and $w_i\, f(w)$ and integrate on $\R^3$ in the $w$-variable.
 We get the identities, when $i\neq j$,
\begin{equation}\label{prin1lfd}  
 v_i\, \frac{\pa_j f}{f \,(1-\dd f)}(v) - v_j \,  \frac{\pa_i f}{f \,(1-\dd f)}(v) = \int_{\R^3}  f(w)\, r_{i,j}^f(v,w) \, dw, 
\end{equation}
\begin{equation}\label{prin2lfd} 
  I_i(f)\,  \frac{\pa_j f}{f \,(1-\dd f)}(v) - K\, v_j + \frac1{\var} \int_{\R^{3}}\ln(1-\dd f(w)) \,w_j\,d w  =  -  \int_{\R^3} w_i\, f(w)\, r_{i,j}^f(v,w) \, dw . 
\end{equation}
This last identity can be rewritten as 
$$   I_i(f)\, \bigg(  \frac{\pa_j f}{f \,(1-\dd f)}(v) - K\, v_j  \bigg) = - (I_i(f) - 1)\, K\, v_j $$
\begin{equation}\label{rew_lfd} 
- \frac1{\var} \int_{\R^{3}}\ln(1-\dd f(w)) \,w_j\,d w   -  \int_{\R^3} w_i\, f(w)\, r_{i,j}^f(v,w) \, dw . 
\end{equation}
Multiplying identity (\ref{prin1lfd}) by $v_i\,v_j\,f(v)\, (1 - \var f(v))$  and integrating on $\R^3$ in the $v$-variable,
we see that
$$ I_j(f) - I_i(f) =  \int\int_{\R^3 \times \R^3}  f(v)\,(1-\dd f(v))\, f(w)\, v_i\, v_j \, r_{i,j}^f(v,w) \, dv dw $$
so that, using again the Cauchy-Schwarz inequality as in the case of the Landau equation, 
$$ | I_j(f) - I_i(f)|
 \le  \sqrt{\frac23}\, \kappa_0^{-1/2}\, D_{\g, LFD}(f)^{1/2} \, s_{\g}(f)^{1/2}, $$
 where $s_{\g}(f)$ is defined by (\ref{vsf}).
As a consequence, once again as in the case of the Landau equation, 
$$ |I_i(f) - 1|^2 
\le \frac8{27}\, \kappa_0^{-1}\,  D_{\g, LFD}(f) \, s_{\g}(f). $$
We deduce from this estimate that
$ I_i(f) \ge \frac12$ as soon as 
\begin{equation}\label{conddplfd}
  D_{\g, LFD}(f) \, s_{\g}(f)  \, \kappa_0^{-1} \le \frac{27}{32} .
\end{equation}
Then, multiplying (\ref{prin2lfd}) by $f(v)$ and integrating on $\R^3$ with respect to the variable $v$, we see that
$$ \frac1{\var} \int_{\R^{3}}\ln(1-\dd f(w)) \,w_j\,d w =  \int\int_{\R^3 \times \R^3}  f(v)\, f(w)\, w_i \, r_{i,j}^f(v,w) \, dv dw ,$$
so that 
$$ \bigg|\frac1{\var} \int_{\R^{3}}\ln(1-\dd f(w)) \,w_j\,d w \bigg| \le   \bigg( \int\int_{\R^3 \times \R^3}  f(v)\, f(w)\, |v-w|^{\g}
 \,| r_{i,j}^f(v,w) |^2 \, dv dw \bigg)^{1/2}$$
$$ \times \,  \bigg( \int\int_{\R^3 \times \R^3}  f(v)\, f(w)\, |v-w|^{-\g} \,| w_i |^2 \, dv dw \bigg)^{1/2} $$
$$ \le \sqrt{2}\, \kappa_0^{-1}\,  D_{\g , LFD}(f)^{1/2} \, \Sigma_{\g}(f)^{1/2}, $$
where
$$ \Sigma_{\g}(f) := \int\int_{\R^3 \times \R^3}  f(v)\, f(w)\, |v-w|^{-\g} \,| w |^2 \, dv dw . $$

Under condition (\ref{conddplfd}), we see that estimate (\ref{rew_lfd}) leads to
$$  \int_{\R^3} \bigg|  \frac{\pa_j f}{f\,(1 - \var f)} (v)  - K\, v_j  \bigg|^2\, f(v) \, dv
  \le 12 \, (I_i(f) - 1)^2\,K^2\, I_j(f)  + 24\, \kappa_0^{-2}\,  D_{\g, LFD}(f) \,   \Sigma_{\g}(f) $$ 
$$ +\,  12\,  \int_{\R^3}  f(v)  \,\bigg[ \int_{\R^3}  f(w)\, |v-w|^{\g} \, |r_{i,j}^f(v,w)|^2 \, dw \bigg] \,dv\,\, \sigma_{\g}(f) ,$$
so that finally, summing from $j= 1$ to $3$ (and choosing $i=j+1$)
$$  \int_{\R^3} \bigg|  \frac{\nabla f}{f\,(1 - \var\,f)} (v) - K\, v  \bigg|^2\, f(v) \, dv \le  \bigg[\frac{32}3\, \kappa_0^{-1}\,K^2 \, s_{\g}(f)  + 72\, \kappa_0^{-2}\,  \Sigma_{\g}(f) $$
$$  + \,24 \,  \kappa_0^{-2}\,  \sigma_{\g}(f) \bigg]\,  D_{\g, LFD}(f) . $$
At this level, we observe that (working as for establishing (\ref{aa2})),
 \begin{equation}\label{aa4}
\Sigma_{\g}(f)  \le  3 \,\left( 1 + 2 \sqrt{\pi}\, ||f||_{L^2(\R^3)} \right) .
\end{equation}
Observing that 
$$  |K| \le \int_{\R^3} \frac{f(v)}{1 - \var\, f(v)}\, dv \le \kappa_0^{-1}, $$ 
 and remembering (\ref{aa2}), (\ref{aa3}) and (\ref{aa4}), we can
 conclude the proof of Proposition \ref{main2}. 
$\square$

\bigskip

 We end up this section by briefly explaining how Proposition \ref{main2} leads to a result of exponential convergence with an explicit and $\var$-independent (for $\var$ not too large) rate towards equilibrium for the solutions to the (spatially homogeneous) Landau-Fermi-Dirac equations with hard potentials or hard spheres (that is $\g \in ]0,1]$). Since the proof closely follows the computations and estimates of \cite{ABDL}, we do not present all details. 
\medskip

First, observing that $\int_{\R^3} f(v)\, |v|^4\, dv$, $||f||_{L^2(\R^3)}$ and $|| \, |\cdot|^2\, f||_{L^2(\R^3)}$ can be estimated by 
$||f||_{L^2_6(\R^3)}$ (see Remark 1), we see thanks to Proposition \ref{main2} that for some $A_f, B_f >0$ depending only on $||f||_{L^2_6(\R^3)}$ 
and $\kappa_0$, the estimate
\begin{equation}\label{bel}
 \int_{\R^3} \bigg| \frac{\nabla f (v)}{f(v) \,(1-\dd f(v))} - K\, v \bigg|^2 \, f(v)\, dv \le B_f\, D_{\gamma, LFD}(f)
\end{equation}
holds for all normalized $f >0$ such that (\ref{kap}) is satisfied, and such that
$$ A_f\, D_{\gamma, LFD}(f) \le 1 . $$
\par
Estimate (\ref{bel}) can be rewritten (see \cite{ABDL}, subsection 3.1) as 
\begin{equation}\label{bel2}
 \int_{\R^3} \bigg| \frac{\nabla f (v)}{f(v) \,(1-\dd f(v))} + 2 b_{\var}\, v \bigg|^2 \, f(v)\, dv  - 3\,(K + 2b_{\var})^2 \le B_f\, D_{\gamma, LFD}(f), 
\end{equation}
where $b_{\var}$ is defined (for $\var$ not too large) in \cite{ABDL}, as the unique real number such that
$$ \int_{\R^3} \frac{ a_{\var}\, e^{- b_{\var}\, |v|^2} }{ 1 + \var\, a_{\var}\, e^{- b_{\var}\, |v|^2} } \, \left(\begin{array}{c} 1\\ |v|^2 \end{array} \right)\, dv = \left(\begin{array}{c} 1\\ 3 \end{array} \right) . $$
Using the estimate (extracted from \cite{ABDL}, Lemma 3.1) 
$$ |K + 2b_{\var}| \le \frac{2 \var}{\kappa_0^2} \, \max( ||f||_{\infty}, ||{\mathcal{M}}_{\var}||_{\infty}) \, ||f - {\mathcal{M}}_{\var} ||_{L^1(\R^3)}, $$
where 
$$ {\mathcal{M}}_{\var} :=  \frac{ a_{\var}\, e^{- b_{\var}\, |v|^2} }{ 1 + \var\, a_{\var}\, e^{- b_{\var}\, |v|^2} } , $$
we obtain the estimate 
$$
 \int_{\R^3} \bigg| \frac{\nabla f (v)}{f(v) \,(1-\dd f(v))} + 2 b_{\var}\, v \bigg|^2 \, f(v)\, dv  - 12\,  \frac{\var^2}{\kappa_0^4} \, \max( ||f||_{\infty}^2, ||{\mathcal{M}}_{\var}||_{\infty}^2) \, ||f - {\mathcal{M}}_{\var}||_{L^1(\R^3)}^2 $$
\begin{equation}\label{bel3}
 \le B_f\, D_{\gamma, LFD}(f).
\end{equation}
Defining 
\begin{equation}\label{aaa}
 H(f|{\mathcal{M}}_{\var}) := S_{\var}({\mathcal{M}}_{\var}) - S_{\var}(f), 
\end{equation}
where 
$$ S_{\var}(f) := -\frac1{\var} \int_{\R^3} \bigg[\var\,f\, \log(\var f) + (1 - \var f)\, \log(1 - \var f)\bigg]\, dv,$$
and
 using the Csiszar-Kullback inequality for Fermi-Dirac relative entropy (cf. \cite{Lu})
$$  ||f - {\mathcal{M}}_{\var}||_{L^1(\R^3)}^2 \le  2\,  H(f|{\mathcal{M}}_{\var}) , $$
which holds for normalized $f$,
 and the corresponding logarithmic Sobolev-like inequality (cf. \cite{cari})
$$   \int_{\R^3} \bigg| \frac{\nabla f (v)}{f(v) \,(1-\dd f(v))} + 2 b_{\var}\, v \bigg|^2 \, f(v)\, dv \ge  2 b_{\var}\, H(f|{\mathcal{M}}_{\var}) , $$ 
we end up with
\begin{equation}\label{bel4}
 \bigg[\, 2 b_{\var} - 24 \,  \frac{\var^2}{\kappa_0^4} \, \max( ||f||_{\infty}^2, ||{\mathcal{M}}_{\var}||_{\infty}^2) \, \bigg] \, H(f|{\mathcal{M}}_{\var}) \le B_f\, D_{\gamma, LFD}(f).
\end{equation}
This last estimate is a functional estimate which holds for all normalized $f>0$ such that (\ref{kap}) is satisfied, under the condition $ A_f\, D_{\gamma, LFD}(f) \le 1$. 
\medskip

We now consider a normalized initial datum  $f_{in} \ge 0$  in $L^1_4(\R^3)$, and a solution $f:= f(t,v)\ge 0$ of the (spatially homogeneous) Landau-Fermi-Dirac $\pa_t f = Q_{\g, LFD}(f,f)$ (for $\g \in ]0,1]$), known to satisfy the estimate 
$$ \sup_{t \ge t_0} ||f(t,\cdot)||_{L^2_6(\R^3)} +  \sup_{t \ge t_0} ||f(t,\cdot)||_{L^{\infty}(\R^3)} \le C_{t_0}, $$
for some $C_{t_0} > 0$ which depends on $t_0 >0$ (cf. Propositions 4.1 and 4.2 of \cite{ABDL}).
\medskip

This ensures that for $\var \in ]0, \var^{*}]$, where $\var^{*} >0$ is given (and can be made explicit), the quantities
$A_{f(t,\cdot)}$, $B_{f(t,\cdot)}$, $||f(t,\cdot)||_{\infty}$ are bounded when $t \ge t_0$, and (\ref{kap}) holds
for $f(t, \cdot)$. 
\medskip

At this level, it is possible to use Lemma \ref{l1} since the entropy identity 
$$ - \frac{d}{dt} H(f(t, \cdot)|{\mathcal{M}}_{\var}) =  D_{\gamma, LFD}(f(t, \cdot)) $$
 holds for the considered solutions
of the  Landau-Fermi-Dirac equation, and to conclude. 
\bigskip

{\bf{Acknowledgement}}: This paper is dedicated to the memory of Maria Conceicao Carvalho, who played a pionneering role in the 
establishment of the entropy-entropy dissipation method in kinetic theory.

\end{document}